\documentclass [11pt]{article}
\usepackage{amsthm,amssymb,amscd}
\theoremstyle{plain}
\newtheorem{rem}{Remark}[section]

\newtheorem{thm}{Theorem}[section]
\theoremstyle{definition}

\newcommand{\bb}[1]{\mbox{$\mathbb{#1}$}}
\begin{document}
\bibliographystyle{plain}

\begin{center}
{\huge\bf A short proof of a formula of Brasselet, L\^{e} and Seade}\\
{\huge\bf for the Euler obstruction}
$ $\\ [2ex]
J\"{o}rg Sch\"{u}rmann\\
Westf. Wilhelms-Universit\"{a}t, SFB 478\\
"Geometrische Strukturen in der Mathematik"\\
Hittorfstr. 27, 48149 M\"{u}nster, Germany\\
E-mail: jschuerm@math.uni-muenster.de
\end{center}
$ $\\ [2ex]
Using the index theorem of Dubson, L\^{e}, Ginsburg and Sabbah
for the vanishing cycle functor (which is a special case of \cite{Sch2}),
we give a short proof of a generalization of the following theorem
\cite[thm.3.1]{BLS}:
\begin{thm} \label{thm:BLS}
Let $(X,0)$ be a germ of an equidimensional complex space (of positive
dimension) in $\bb{C}^{N}$. Let $V_{i}, \;i=1,\cdots,k$, be the (connected)
strata of a Whitney stratification of a small representative $X$ of $(X,0)$.
Then there is an open dense Zariski subset $\Omega$ in the space of complex
linear forms $\bb{L}$ on $\bb{C}^{N}$, such that for every $l\in \Omega$, 
there is 
$\epsilon_{0}$, such that for any $\epsilon$, $\epsilon_{0}> \epsilon >0$
and $t_{0}$ sufficiently small, we have the following formula for the
Euler obstruction of $(X,0)$:
\[Eu(X,0) = \sum_{1}^{k} \chi\bigl(V_{i}\cap \bb{B}_{\epsilon}\cap 
l^{-1}(t_{0})\bigr) \cdot Eu_{V_{i}}(X) ,\]
where $\chi$ denotes the Euler-Poincar\'{e} characteristic and $Eu_{V_{i}}(X)$
is the value of the Euler obstruction of $X$ at any point of $V_{i},\; 
i=1,\cdots,k$.
\end{thm}

Here $\bb{B}_{\epsilon}$ is the corresponding closed ball in $\bb{C}^{N}$,
and $\chi$ is the usual Euler characteristic 
\[\chi(Y):=\chi\bigl(H^{*}(Y,\bb{Q})\bigr)\]
of a topological space $Y$ with finite dimensional cohomology 
$H^{*}(Y,\bb{Q})$. Similarly, let us denote by $\chi_{c}(Y)$ the Euler 
characteristic of $H^{*}_{c}(Y,\bb{Q})$ for locally compact $Y$ (if defined).

\begin{rem} This result goes already back to an old preprint \cite{Du}
of Dubson (compare with \cite[thm.V.1.1,p.69]{Du} (R\'{e}duction \`{a}
la section hyperplane) and \cite[thm.IV.4.1.2,p.67]{Du} (Lemme de Hopf 
singulier)), which was not mentioned in \cite{BLS}.
\end{rem} 

Consider a Whitney stratification of a germ of a complex space $(X,0)$
of positive dimension in $\bb{C}^{N}$ as above, but this time whe don't 
assume $(X,0)$ is equidimensional. Let us refine the stratification 
by adding the point stratum $V_{0}:=\{0\}$, and then assume $0\not\in V_{i}$
for $i>0$. Take a linear form
$l\in \bb{L}$, such that $dl_{0}$ is a nondegenerate covector in the sense
of stratified Morse theory (i.e. it does not vanish on any generalized tangent 
space
$lim_{x_{n}\to 0} T_{x_{n}}V_{i}$ at $0$, with $x_{n}$ a sequence in $V_{i}$,
$i=1,\cdots,k$). This condition defines a Zariski open dense subset $\Omega$
of $\bb{L}$. Then 
\[V(i,\epsilon,t_{0}):= V_{i}\cap \bb{B}_{\epsilon}\cap l^{-1}(t_{0})\]
is for $i>0$ a smooth manifold with boundary for $0<|t_{0}|<<\epsilon_{0}<<1$
(this is also true for any $l\in \bb{L}$, since one can refine our 
stratification to an $a_{f}$-stratification).
By Poincar\'{e}-duality we get
\begin{equation}
\chi(V(i,\epsilon,t_{0}))
= (-1)^{2(dim_{C}(V_{i})-1)}\cdot \chi_{c}(V^{\circ}(i,\epsilon,t_{0}))
= \chi_{c}(V^{\circ}(i,\epsilon,t_{0})) ,
\end{equation}
with $V^{\circ}(i,\epsilon,t_{0}):=V_{i}\cap \bb{B}^{\circ}_{\epsilon}\cap 
l^{-1}(t_{0})$, and $\bb{B}^{\circ}_{\epsilon}$ the corresponding open ball.\\

Consider now a function $\alpha: X\to \bb{Z}$, which is constructible with
respect to the refined Whitney stratification (i.e. is constant on the 
strata). Let $i_{0}: \{0\}\to X\cap\{l=0\}$ and $i: X\cap\{l=0\}\to X$ be
the inclusions. 
Then we have for $0<|t_{0}|<<\epsilon_{0}<<1$:
\begin{equation} \label{eq:costalk}
i_{0}^{!}\bigl(\psi_{l}(\alpha)\bigr)(0) = 
\sum_{1}^{k} \chi_{c}(V^{\circ}(i,\epsilon,t_{0}))\cdot \alpha_{i} ,
\end{equation}
with $\alpha_{i}$ the value of $\alpha$ on $V_{i}$, and $\psi_{l}$
the nearby cycle functor with respect to $l$ on the level of constructible
functions (and similarly for $i_{0}^{!}$). This can be defined by two different
(but equivalent) approaches:

\begin{itemize}
\item One can use the  sophisticated language of analytically constructible
complexes of sheaves:\\   
Choose a bounded complex ${\cal F}$ of sheaves of $\bb{C}$-vector spaces on 
$X$, which is constructible with respect to our refined Whitney 
stratification, such that all stalk complexes ${\cal F}_{x}$ have finite 
dimensional cohomology with $\chi\bigl(H^{*}({\cal F}_{x})\bigr)
= \alpha(x)$ for all $x\in X$. Then
\begin{equation}
i_{0}^{!}\bigl(\psi_{l}(\alpha)\bigr)(0):=
\chi\bigl(H^{*}(i_{0}^{!}\psi_{l}({\cal F}))\bigr) .
\end{equation}
This is well defined by \cite[Part II.3, Part VI]{Sch},
and equation (\ref{eq:costalk}) follows from \cite[Part V, lem.3.4]{Sch}.
\item By linearity one can assume $\alpha=1_{Z}$, with $Z$ a closed union
of strata and then one uses a direct geometric definition in terms of a
corresponding local Milnor fibration:
\begin{equation} \label{eq:costalk2}
i_{0}^{!}\bigl(\psi_{l}(1_{Z})\bigr)(0) = 
\chi_{c}(M^{\circ}(l|Z)) =
\sum_{V_{i}\subset Z} \chi_{c}(V^{\circ}(i,\epsilon,t_{0})) ,
\end{equation}
with $M^{\circ}(l|Z) := Z\cap \bb{B}^{\circ}_{\epsilon}\cap 
l^{-1}(t_{0})$ 
(for  $0<|t_{0}|<<\epsilon_{0}<<1$) the open local Milnor fiber of $l|Z$
at $0$.
\end{itemize}

But by \cite[Part VI, ex.2.3]{Sch} we have
\begin{equation}
i_{0}^{!}\bigl(\psi_{l}(\alpha)\bigr)(0) =
i_{0}^{*}\bigl(\psi_{l}(\alpha)\bigr)(0) .
\end{equation}
In our situation this formula can also be explained more geometrically
in the following way: Choose a Whitney stratification adapted to $\alpha$
and $\{l=0\}$ (e.g. $\{l=0\}$ is a union of strata), which is also an
$a_{f}$-stratification. Then the corresponding local Milnor fibrations
are "locally constant" along the strata of $\{l=0\}$ so that
$\psi_{l}(\alpha)$ is constructible with respect to the induced stratification
of $\{l=0\}$. So by linearity it is again enough to show
\begin{equation}
i_{0}^{!}\bigl(1_{Z}\bigr)(0) =
i_{0}^{*}\bigl(1_{Z}\bigr)(0) 
\end{equation}
for $Z$ a closed union of strata of $\{l=0\}$. But 
$i_{0}^{*}\bigl(1_{Z}\bigr)(0) -
i_{0}^{!}\bigl(1_{Z}\bigr)(0)$ is 
by the locally conic structure of $Z$ just the Euler characteristic of
the link $Z\cap\{|z|=\epsilon\}$ of $Z$ in $0$ (for $0<\epsilon<<1$),
which vanishes by a classical result of Sullivan (compare \cite{BV,Su}
and also  \cite[Ex.IX.12, p.409]{KS}).\\

From the definition of the vanishing cycle functor, we get the equality
\begin{equation}
\alpha(0) = i_{0}^{*}i^{*}(\alpha)(0) =
i_{0}^{*}\bigl(\psi_{l}(\alpha)\bigr)(0) ,
\end{equation}
if $\chi\bigl(H^{*}(i_{0}^{*}\phi_{l}({\cal F}))\bigr) = 0$.
Or more directly in terms of constructible functions:
\[ i_{0}^{*}\bigl(\phi_{l}(\alpha)\bigr)(0) :=
i_{0}^{*}\bigl(\psi_{l}(\alpha)\bigr)(0) - i_{0}^{*}i^{*}(\alpha)(0)\;.\]

But this follows from the index theorem of Dubson, L\^{e}, Ginsburg and Sabbah
\cite{Du2,Le,Le2,Gi,Sa} for the vanishing cycle functor, if the graph of $dl$ 
doesn't intersect
the support $|CC({\cal F})|=|CC(\alpha)|$ of the characteristic cycle of 
${\cal F}$ or $\alpha$ in $T^{*}\bb{C}^{N}$ over an open neighborhood $U$ of 
$0$. More precisely, this index formula implies (\cite[thm.4.5,rem.4.6]{Sa}):
\begin{equation}
-i_{0}^{*}\bigl(\phi_{l}(\alpha)\bigr)(0) =
-\chi\bigl(H^{*}(i_{0}^{*}\phi_{l}({\cal F}))\bigr) =
\sharp\bigl(dl(U)\cap CC({\cal F})\bigr) , 
\end{equation}
if $dl(U)\cap |CC({\cal F})|\subset \{(0,dl_{0})\}$. \\

\begin{rem} For another purely real proof of this intersection formula
compare with \cite{Sch2} (which is also closely related to 
\cite[thm.9.5.6]{KS}).
\end{rem} 
  
Let $Eu(cl(V_{i}),\cdot)$ be the Euler obstruction of the closure of the
stratum $V_{i}$ ($i>0$). This function is constructible with respect to
our refined Whitney stratification (compare \cite[part VI]{Sch}
or \cite[exp.6, cor.10.2, p.125]{CEP}), with 
\[|CC(Eu\bigl(cl(V_{i}),\cdot)\bigr)| = cl\bigl(T^{*}_{V_{i}}\bb{C}^{N}\bigr) .\]
But the condition $(0,dl_{0})\not\in cl\bigl(T^{*}_{V_{i}}\bb{C}^{N}\bigr)$ is
equivalent to the fact, that $dl_{0}$ 
does not vanish on any generalized tangent space
$lim_{x_{n}\to 0} T_{x_{n}}V_{i}$ at $0$, with $x_{n}$ a sequence in $V_{i}$
($i\in \{1,\cdots,k\}$). \\

Altogether we get the following generalization of theorem 
\ref{thm:BLS} of Brasselet, L\^{e} and Seade:
\begin{thm}
Let $(X,0)$ be a germ of a complex space (of positive
dimension) in $\bb{C}^{N}$. Let $V_{i}, \;i=1,\cdots,k$, be the (connected)
strata of a Whitney stratification of a small representative $X$ of $(X,0)$.
Refine the stratification 
by adding the point stratum $V_{0}:=\{0\}$, and then assume $0\not\in V_{i}$
for $i>0$. Fix a (constructble) function $\alpha$, which is a linear 
combination of $Eu(cl(V_{i}),\cdot)$ (for some $i\in I\subset \{1,\dots,k\}$).
Consider a complex linearform $l\in \bb{L}$, with $dl_{0}$ 
not vanishing on any generalized tangent space
$lim_{x_{n}\to 0} T_{x_{n}}V_{i}$ at $0$, with $x_{n}$ a sequence in $V_{i}$
($i\in I$). Then there is an $\epsilon_{0}>0$ with
\[\alpha(0) = \sum_{1}^{k} \chi_{c}\bigl(V^{\circ}(i,\epsilon,t_{0})\bigr)
\cdot \alpha_{i} 
= \sum_{1}^{k} \chi\bigl(V_{i}\cap \bb{B}_{\epsilon}\cap l^{-1}(t_{0})\bigr)
\cdot \alpha_{i} ,\]
for $0<|t_{0}|<<\epsilon<\epsilon_{0}$.
\end{thm}
 
\begin{rem} The same result is true with the same proof, if the
stratification satisfies suitable weaker regularity conditions than
Whitney b-regularity, e.g. the stratification is (locally) C-regular  
in the sense of Bekka \cite{Be} (compare with \cite[Part VI]{Sch}, even for
a more general class of stratifications). Moreover, this formula is more
generally true (with the same proof) for any holomorphic function germ
$l: (\bb{C}^{N},0)\to (\bb{C},0)$, with $dl_{0}$ as before.
\end{rem}

\end{document}